\documentclass[12pt]{amsart}
\usepackage{amsfonts}
\usepackage{amssymb}
\theoremstyle{definition}

\theoremstyle{remark}

\newcommand{\E}{\mathbb E}

\begin{document}

\centerline{Riv. Mat. Univ. Parma \ (5) \ {\bf 5} (1996), 15-21}

\vspace{0.8in}

\centerline{\large\bf A characterization of the extrinsic spheres}
\centerline{\large\bf in a Riemannian manifold
\footnote{Received November 28, 1995. AMS classification 53 C 40. The research has been
partially supported by the Ministry of Education of Bulgaria, Contract MM 313/94.}}

\vspace{0.8in}
\centerline{\large Ognian  Kassabov
\footnote{Higher Transport School "T. Kableschkov", Section of Mathematics, Slatina, 
1574 Sofia, Bulgaria.}}

\vspace{1.in}
{\bf 1 - Introduction} 

\vspace{0.3in}
In [{\bf 3}] K. Ogiue and R. Takagi propose different criteria for a surface in $\E^3$
to be a sphere and give the following useful and practical condition

\vspace{0.3in}
T\,h\,e\,o\,r\,e\,m A. {\it Let $M$ be a surface in $\E^3$. Suppose that, through each
point $p\in M$ there exist two circles of $\E^3$ such that:

\ \ \ they are contained in $M$ in a neighbourhood of $p$

\ \ \ they are tangent to each other at $p$.

Then $M$ is locally a plane or a sphere.}

\vspace{0.3in}
The so-called extrinsic spheres are a natural generalization of the ordinary sphere
in $\E^m$. It is interesting to have a characterization of an $n$-sphere, or more
generally of an extrinsic sphere, similar to Theorem A. K. Ogiue and R. Takagi
give such a criterion in [{\bf 3}] by means of $n^2$ circles through each point of the
submanifold. Here we propose an analogue of Theorem A for extrinsic spheres. Namely we prove

\vspace{0.3in}
T\,h\,e\,o\,r\,e\,m B. {\it Let $M$ be an $n$-dimensional $(n>2)$ submanifold of a
Riemannian manifold $\widetilde M$. Then $M$ is either a totally geodesic submanifold of 
$\widetilde M$ or an extrinsic sphere of $\widetilde M$, if through each point $p$ of $M$ there 
exist two $(n-1)$-dimensional extrinsic spheres of $\widetilde M$, such that: 

\ \ \  they are contained in $M$ in a neighbourhood of $p$

\ \ \  they are tangent to each other at $p$.}

\vspace{0.3in}
C\,o\,r\,o\,l\,l\,a\,r\,y. {\it Let $M$ be an $n$-dimensional $(n>2)$ submanifold of 
the Euclidean space $\E^m$. Suppose that through each point $p$ of $M$ there exist
two $(n-1)$-spheres of $\E^m$, such that:

\ \ \ they are contained in $M$ in a neighbourhood of $p$

\ \ \ they are tangent to each other at $p$.

Then $M$ is locally an $n$-plane or an $n$-sphere in $\E^m$.}

\vspace{0.8in}
{\bf 2 - Preliminaries}

\vspace{0.3in}

Let $\widetilde M$ be a Riemannian manifold with metric tensor $g$ and let $M$ be an
$n$-dimensional submanifold of $\widetilde M$. Denote by $\widetilde\nabla$ and $\nabla$ 
the Riemannian connections of $\widetilde M$ and $M$, respectively. Then the {\it Gauss
formula} is

$$
	\widetilde\nabla_XY = \nabla_XY +\sigma(X,Y)
$$

\noindent
for all vector fields $X,\,Y$ on $M$, where $\sigma$ is the {\it second fundamental
form} of $M$ in $\widetilde M$. Let $\xi$ be a normal vector field. Then the {\it Weingarten
formula} is

$$
	\widetilde\nabla_X\xi = -A_{\xi}X + D_X\xi
$$

\noindent
where $-A_\xi X$ and $D_X\xi$ are the tangential and the normal components of $\widetilde\nabla_X\xi$
respecti\-vely. Usually $A_\xi$ is called the {\it shape operator}, corresponding to $\xi$ and
D the {\it connection in the normal bundle}. Also $g(A_\xi X,Y)=g(\sigma(X,Y),\xi)$ holds
good. Recall then that the covariant derivative of $\sigma$ with respect to the {\it connection
$\overline\nabla$ of van der Waerden-Bortolotti} is given by

$$
	(\overline\nabla_X\sigma)(Y,Z)=D_X\sigma(Y,Z)-\sigma(\nabla_XY,Z)-\sigma(Y,\nabla_XZ) \,.
$$  

The {\it mean curvature vector} $H$ of $M$ in $\widetilde M$ is defined by $H=n^{-1}\,{\rm trace}\,\sigma$.

The submanifold $M$ is called {\it totally umbilical}, if $\sigma(x,y)=g(x,y)H$
for all $x,\,y \in T_pM$, $p\in M$ or equivalently $A_\xi x=g(\xi,H)x$ for all
$x\in T_pM$, $\xi\in (T_pM)^{\perp}$, $p\in M$. In particular, if $\sigma$ vanishes
identically, $M$ is said to be a {\it totally geodesic} submanifold of $\widetilde M$.
A normal vector field $\xi$ is called {\it parallel}, if $D_X\xi=0$ for any vector
field $X$ on $M$. 

A regular curve $\tau=(x_s)$ parametrized by arc length $s$ is called a {\it circle} in
$\widetilde M$, if there exists a field $Y_s$ of unit vectors along $\tau$ and a 
positive constant $k$, such that

$$
	\widetilde\nabla_{X_s}X_s=kY_s\,, \qquad  \widetilde\nabla_{X_s}Y_s=-kX_s\,,
$$

\noindent
where $X_s$ denotes the tangent vector of $\tau$ [{\bf 2}]. The number $k$ is called the 
{\it radius} and $Y_s$ the {\it main normal} of $\tau$.

A submanifold $M$ of $\widetilde M$ is said to be an {\it extrinsic sphere}, if it 
is totally umbilical and has non-zero parallel mean curvature vector. For an extrinsic 
sphere $M$ we have

$$
	\widetilde\nabla_X\widetilde\nabla_XX=\nabla_x\nabla_XX-g(H,H)X
$$

\noindent
and hence

$$
	g(\widetilde\nabla_X\widetilde\nabla_XX,\xi)=0   \leqno (2.1)
$$

\noindent
for any unit vector field $X$ on $M$ and any vector field $\xi$ normal to $M$.

\vspace{0.8in}
{\bf 3 - General lemmas}

\vspace{0.3in}
In this section we prepare two lemmas (Lemma 2 and Lemma 3) which will be useful
in the proof of our theorem. The first of them shows that an extrinsic sphere
is determined (locally) by its second fundamental form at one point. The second
gives a condition for the second fundamental forms of two hypersurfaces to coincide.

\vspace{0.3in}
L\,e\,m\,m\,a 1. {\it Let $S$ be an extrinsic sphere in a Riemannian manifold $\widetilde M$
and denote by $H$ its mean curvature vector. Then every geodesic of $S$ through a point
$p$ of $S$ is a circle in $\widetilde M$ of radius $k=(g(H,H))^{\frac12}$ and its mean normal
vector at $p$ is $k^{-1}H_p$.}

\vspace{0.3in}
P\,r\,o\,o\,f. Note that $k=(g(H,H))^{\frac12}$ is a non-zero constant since $H$ is parallel
and non-zero. Let $\tau=(x_s)$ be a geodesic of $S$ parametrized by arc length. Denote by 
$X_s$ its tangent vector and define $Y_s$ by $H_{x(s)}=kY_s$. Then by the Gauss and
Weingarten formulas

$$
	\widetilde\nabla_{X_s}X_s=kY_s   \qquad	\qquad  \widetilde\nabla_{X_s}Y_s=-kX_s 
$$

hold good and the lemma is proved.

\vspace{0.3in}
L\,e\,m\,m\,a 2. {\it Let $\widetilde M$ be a Riemannian manifold. Suppose that through 
a point $p$ of $\widetilde M$ there exist two extrinsic spheres $S_1$ and $S_2$ of 
$\widetilde M$, which are tangent to each other at $p$. If their second fundamental forms 
at $p$ coincide, then $S_1$ and $S_2$ coincide in a neighbourhood of $p$.}

\vspace{0.2in}
P\,r\,o\,o\,f. Denote by $H_i$ the mean curvature vector of $S_i$ in $\widetilde M$. Note
that $H_1=H_2$ at $p$. Let $r$ be a positive number, such that ${\rm exp}_i$ is a diffeomorphism
of a neighbourhood $N_i(p,r)$ of the origin of $T_pS_1=T_pS_2$ and a neighbourhood 
$U_i(p,r)$ of $p$ in $S_i$, $i=1,2$, as in Proposition 3.4, Chapter IV of [{\bf 1}].

For a point $q \in U_1(p,r)$ let $\tau=(x_s)$, $s \in [0,s_0]$ be the only geodesic of $S_1$
in $U_1(p,r)$ joining $p$ and $q$ and parametrized by arc length. Let $\bar\tau=(\bar x_s)$ 
be defined by $\bar x_s={\rm exp}_2({\rm exp}_1^{-1}x_s)$. Then $\bar\tau $ is a geodesic in 
$S_2$ through $p$.

According to Lemma 1, $\tau$ and $\bar\tau$ are circles in $\widetilde M$ through $p$ of radius
$k=(g(H_1,H_1))^{\frac12}$ and its main normal at $p$ is $k^{-1}(H_1)_p=k^{-1}(H_2)_p$. But
such a circle is unique [{\bf 2}]. Then $\tau$ and $\bar\tau$ coincide locally.

Let $\bar q={\rm exp}_2({\rm exp}_1^{-1}q)$. Since $d(p,q)=s_0=d(p,\bar q)$ it follows
$q=\bar q$. So $q \in S_2$, thus proving the lemma. 

\vspace{0.3in}
L\,e\,m\,m\,a 3. {\it Let $S_1$ and $S_2$ be two hypersurfaces through a point $p$ of
a Riemannian manifold $M$. Suppose also that $S_1$ and $S_2$ are tangent  to each other
and that there exists a smooth unit vector field $N$, defined in a neighborhood of
$p$ in $M$ and such that $N$, resticted to $S_i$, is normal to $S_i$, for $i=1,\,2$.
Then the second fundamental forms of $S_1$ and $S_2$ coincide at $p$.}

\vspace{0.2in}
P\,r\,o\,o\,f. Let $x \in T_pS_1=T_pS_2$. Since $S_1$ and $S_2$ are hypersurfaces the
Weingarten formulas for $S_1$ and $S_2$ in $M$ imply $\nabla_xN|_{S_i}=-A^i_{N_p}x $ for
$i=1,\,2$, where $A^i_{N_p}$ is the shape operator of $S_i$ in $M$ and $N|_{S_i}$
denotes the restriction of $N$ to $S_i$. But $\nabla_xN=\nabla_xN|_{S_1}=\nabla_xN|_{S_2}$.
Hence $A^1_{N_p}=A^2_{N_2}$ and consequently the second fundamental forms of $S_1$
and $S_2$ coincide at the point $p$.

\vspace{0.8in}
{\bf 4 - Proof of Theorem B}

\vspace{0.3in}
For a point $p$ of $M$ let $S_{1p}$ and $S_{2p}$ be extrinsic spheres through $p$ as in 
the statement of the theorem. Let $X$ be a vector field on $S_{1p}$. Then a direct calculation gives

$$
	g(\widetilde\nabla_X\widetilde\nabla_XX,\bar\xi) =
	g((\overline\nabla_X\sigma)(X,X)+3\sigma(\nabla_XX,X),\bar\xi)
$$

\noindent
for any vector field $\bar\xi$, normal to $M$. Because of (2.1) this implies

$$
	(\overline\nabla_X\sigma)(X,X)+3\sigma(\nabla_XX,X)=0   \leqno (4.1)
$$

\noindent
for any unit vector field $X$ on $S_{1p}$.

Denote by $\nabla^1$ the Riemannian connection of $S_{1p}$. Let $\xi_1$ be 
a local normal unit vector field for $S_{1p}$ in $M$. Then for all vector fields
$X,\,Y$ on $S_{1p}$ we have

$$
	\nabla_XY=\nabla^1_XY+h_1(X,Y)\xi_1 
$$

\noindent
$h_1\xi_1$ being the second fundamental form of $S_{1p}$ in $M$. Hence, using the
Gauss formula of $M$ in $\widetilde M$ we find 

$$
	\widetilde\nabla_XY=\nabla^1_XY+h_1(X,Y)\xi_1+\sigma(X,Y) \,. 
$$

\noindent
Denote by $H_1$ the mean curvature vector of $S_{1p}$ in $\widetilde M$. Since $S_{1p}$
is totally umbilical in $\widetilde M$, we get

$$
	g(X,Y)H_1=h_1(X,Y)\xi_1+\sigma(X,Y) \ .  \leqno (4.2)
$$

\noindent
Note that $\xi_1$ is orthogonal to $\sigma(X,Y)$ for all vector fields $X,\,Y$ on $S_{1p}$.
Then (4.2) shows that $S_{1p}$ is totally umbilical in $M$. More explicitly, putting 
$\lambda_1=g(H_1,\xi_1)$, from (4.2) we obtain

$$
	h_1(X,Y)=\lambda_1g(X,Y)             \leqno (4.3)
$$

$$
	\sigma(X,Y)=(H_1-\lambda_1\xi_1)g(X,Y)  \,.   \leqno (4.4)
$$

Suppose now that $X$ is a unit vector field on $S_{1p}$. Then $X$ is orthogonal to
$\nabla^1_XX$ and consequently from (4.1), (4.3) and (4.4) we derive

$$
	(\overline\nabla_X\sigma)(X,X)+3\lambda_1\sigma(X,\xi_1)=0   \,.    \leqno (4.5)
$$

Put $\xi=\xi_1(p)$. Then by (4.5) we get

$$
	(\overline\nabla_x\sigma)(x,x)+3\lambda_1(p)\sigma(x,\xi)=0 
$$

\noindent
for any unit vector $x$ in $T_pS_{1p}=T_pS_{2p}$. It is easy to see that we 
have also

$$
	(\overline\nabla_x\sigma)(x,x)+3\lambda_2(p)\sigma(x,\xi)=0 
$$

\noindent
where $\lambda_2(p)\xi g$ is the second fundamental form at $p$ of $S_{2p}$ in $M$.
The last two equations and Lemma 2 imply
$$
	\sigma(x,\xi)=0          \leqno (4.6)
$$ 

\noindent
for any vector $x$ in $T_pS_{1p}=T_pS_{2p}$.

According to (4.4) we may write \
$\sigma(x,y) = g(x,y)\eta$ \ for \ $x,\,y \in T_pS_{1p}\,$, 
where $\eta$ is given by

$$
	\eta = H_1(p) - \lambda_1(p)\xi=H_2(p)-\lambda_2(p)\xi
$$

\noindent
$H_2$ being the mean curvature vector of $S_{2p}$ in $\widetilde M$. Denote
$\zeta=\sigma(\xi,\xi)$. Then using (4.6) we find for any $x,\,y \in T_pM$

$$
	\sigma(x,y)=(g(x,y)-g(x,\xi)g(y,\xi))\eta+g(x,\xi)g(y,\xi)\zeta
$$

\noindent
and hence

$$
	A_{{\rm tr}\,\sigma}x=g(\eta,{\rm tr}\,\sigma)x+g(\zeta-\eta,{\rm tr}\,\sigma)g(x,\xi)\xi  \leqno(4.7)
$$

\noindent
for any $x\in T_pM$. We put \
$\mu_0=g(\eta,{\rm tr}\,\sigma),\,\nu_0=g(\zeta,{\rm tr}\,\sigma)$.
Then according to (4.7) $\xi$ is an eigenvector of \ $A_{{\rm tr}\,\sigma}$ \ at $p$ 
and the corresponding eigenvalue is $\nu_0$. Analogously any unit vector $x$ in \
$T_pS_{1p}=T_pS_{2p}$ \ is an eigenvector of \ $A_{{\rm tr}\,\sigma}$ \ with corresponding
eigen\-value $\mu_0$.

Suppose that $M$ is not totally umbilical at $p$, i.e. $\mu_0\ne \nu_0$. Let the continuous
functions $\mu$ and $\nu$ be eigenvalues of $A_{{\rm tr}\,\sigma}$, such that
$\mu(p)=\mu_0$ and $\nu(p)=\nu_0$. Then $\mu(q)\ne \nu(q)$ for any $q$ in a
sufficiently small neighbourhood $U$ of $p$ in $M$. It follows directly or by the
implicit function theorem that $\nu$ is a smooth function on $U$. Then its corresponding
field $N$ of unit eigenvectors is also smooth on $U$.

We shall show that the restriction of $N$ to $S_{ip}$ is orthogonal to $S_{ip}$ for $i=1,\,2$.
Indeed according to the definitions of $\mu,\,\nu$ and $N$

$$
	A_{{\rm tr}\,\sigma}x =\mu(q)x + (\nu(q)-\mu(q))g(x,N)N_q
$$  

\noindent
for all $x \in T_qM$, $q \in U$. Hence

$$
	g(\sigma(x,y),{\rm tr}\,\sigma)=(\nu(q)-\mu(q))g(x,N)g(y,N)
$$

\noindent
for all orthogonal vectors $x,\,y \in T_qM$, $q\in U$. On the other hand by (4.4) we have
$g(\sigma(x,y),{\rm tr}\,\sigma)=0$ for all orthogonal vectors $x,\,y \in T_qS_{1p}$, $q\in S_{1p}$.
Consequently, for arbitrary orthogonal vectors  $x,\,y \in T_qS_{1p}$, 
$g(x,N)g(y,N)=0$ holds good, i.e. at least one of any two orthogonal vectors in $T_qS_{1p}$ is 
orthogonal to $N_q$. Hence it follows easily that in fact any vector in $T_qS_{1p}$ is
orthogonal to $N_p$. So the restriction of $N$ to $S_{1p}$ is orthogonal to $S_{1p}$. 
Analogously the restriction of $N$ to $S_{2p}$ is orthogonal to $S_{2p}$.

Then, according to Lemma 3 the second fundamental forms of $S_{1p}$ and $S_{2p}$ coincide
at $p$. Now using Lemma 2 we conclude that $S_{1p}$ and $S_{2p}$ coincide in a neighbourhood
of $p$, which is a contradiction. So $\mu_0=\nu_0$ and $M$ is totally umbilical at $p$.
Since $p$ is an arbitrary point of $M$, it follows that $M$ is totally umbilical in
$\widetilde M$.

It remains to show that the mean curvature vector $H$ of $M$ is parallel. Since $M$ is 
totally umbilical, by (4.1) we conclude that $D_xH=0$ for any $x\in T_qS_{1p}$, $q\in S_{1p}$.
Suppose that for the above defined vector $\xi$ we have $D_\xi H\ne 0$. Let $\bar\eta$ be
a normal vector field of $M$, defined in a neighbourhood of $p$, such that the differential
form $\omega$, given by $\omega(X)=g(D_XH,\bar\eta)$ for a vector field $X$ on $M$ is not
zero at any point of $U$. Define a vector field $Y$ on $U$ by $\omega(X)=g(X,Y)$. Then $Y$
does not vanish in $U$. Define $N=(g(Y,Y))^{-\frac12}Y$. Note that the restriction of $N$
to $S_{ip}$ is normal to $S_{ip}$ for $i=1,\,2$. Using again Lemmas 2 and 3 we conclude that 
$S_{1p}$ and $S_{2p}$ coincide in a neighbourhood of $p$, which is a contradiction. So $H$
is parallel. Then $M$ is totally geodesic or an extrinsic sphere in $\widetilde M$,
according to the length of $H$ being zero or not. This proves our theorem.

\vspace{0.8in}
\centerline{\bf References}

\vspace{0.2in}
\noindent
[{\bf 1}] \ \ \ S. K\textsc{obayashi} and K. N\textsc{omizu}, {\it Foundations of differential geometry}, {\bf 1}, 
Interscience, 

\hspace{0.1in} New York 1963.

\noindent
[{\bf 2}] \ \ \ K. N\textsc{omizu} and K. Y\textsc{ano}, {\it On circles and spheres in Riemannian geometry}, Math. 

\hspace{0.1in} Ann. {\bf 210} (1974), 163-170.

\noindent
[{\bf 3}] \ \ \ K. O\textsc{giue} and R. T\textsc{akagi}, {\it A submanifold which contains many extrinsic circles},

\hspace{0.1in}
Tsukuba J. Math. {\bf 8} (1984), 171-182.

\vspace{0.8in}
\centerline{S o m a r i o}

\vspace{0.2in}

{\it Nel 1984 K. Ogiue ed R. Takagi hanno dato una condizione perch\'e una superficie
dello spacio ordinario sia localmente un piano o una sfera. Viene qui ottenuta una 
condizione dello stesso tipo perch\'e una sottovariet\`a $M$ di una variet\`a Riemanniana $ \widetilde M$ sia
totalmente geodesica oppure sia una sfera estrinseca.}

\end{document}